\documentclass[reqno, 11pt]{amsart}
\usepackage{amsfonts}
\usepackage[dvips]{color}
\usepackage{amssymb,latexsym}
\usepackage{a4}
\usepackage{fancybox}
\usepackage{epsfig}
\usepackage[english]{babel}
\usepackage{textcomp}
\usepackage{amsmath}
\usepackage{graphicx}

\newtheorem{theorem}{Theorem}[section]
\newtheorem{proposition}[theorem]{Proposition}
\newtheorem{lemma}[theorem]{Lemma}
\newtheorem{corollary}[theorem]{Corollary}
\newtheorem{definition}[theorem]{Definition}
\newtheorem{remark}{Remark}[section]

\newcommand{\oo}{\infty}

\newcommand{\comp}{\circ}

\makeatletter
\@addtoreset{equation}{section}
\makeatother

\renewcommand\thefigure{\thesection.\@arabic\c@figure}
\renewcommand\thetable{\thesection.\@arabic\c@table}

\newcommand{\1}{\mathbb{1}}

\newcommand{\N}{\mathbb{N}}

\def\1{\,\rlap{\mbox{\small\rm 1}}\kern.15em 1}

\def\build#1_#2^#3{\mathrel{\mathop{\kern 0pt#1}\limits_{#2}^{#3}}}
\def\tend#1#2{\build\hbox to 12mm{\rightarrowfill}_{#1\rightarrow #2}^{ }}

\def\converge#1#2#3#4{\build\hbox to
#1mm{\rightarrowfill}_{#2\rightarrow #3}^{\hbox{\scriptsize #4}}}

\newcommand{\beq}{\begin{equation}}
\newcommand{\eeq}{\end{equation}}

\bibdata{prob}
\bibstyle{alpha}

\title{Dynamical properties of random walks}

\author{Ali Messaoudi$^1$, Glauco Valle$^2$}
\thanks{1. Supported by CNPq grant 307776/2015-8 and  Fapesp project 2013/23643-4}
\thanks{2. Supported by CNPq grant 305805/2015-0}

\address{
\newline
Ali Messaoudi
\newline
UNESP - Departamento de matem\'atica do Instituto de Bioci\^encias Letras e Ci\^encias Exatas.
\newline Rua Crist\'ov\~ao Colombo, 2265, Jardim Nazareth,
15054-000 - S\~ao Jos\'e do Rio Preto, SP, Brasil.
\newline
e-mail: {\rm \texttt{messaoud@ibilce.unesp.br}}
\newline
\newline
Glauco Valle
\newline
Universidade Federal do Rio de Janeiro, Instituto de Matem\'atica.
\newline  Caixa Postal 68530, 21945-970, Rio de Janeiro, Brasil.
\newline
e-mail: {\rm \texttt{glauco.valle@im.ufrj.br}}
}

\subjclass[2010]{primary 47A16 secondary 	60J10}
\keywords{ Random Walks, Discrete Birth and death processes, Linear Dynamics, Supercyclicity, Devaney chaos}

\begin{document}

\maketitle

\begin{abstract}
In this paper, we study  dynamical properties as hypercyclicity,  supercyclicity, frequent hypercyclicity and chaoticity for transition operators associated to countable irreductible Markov chains. As particular cases, we consider simple random walks on $\mathbb{Z}$ and $\mathbb{Z}_+$.
\end{abstract}

\section{Introduction}

Let $X$ be a Banach separable space on $\mathbb{C}$ and $T: X \to X$ be a linear operator on X. The study of the linear dynamical system $(X, T)$ became very active after 1982. Since then related works have built connections between dynamical systems, ergodic theory and functional analysis.  We refer the reader to the
books \cite{bm,em} and to the more recent papers
\cite{BerBonMulPer13,BerBonMulPer14,BesMenPerPui16,SGriEMat14,
SGriMRog14}, where many additional references can be found.

The objective of this paper is to study some central properties of linear dynamical systems as hypercyclicity, supercyclicity, frequent hypercyclicity, and chaoticity among others, for Markov chain transition operators associated  to countable irreductible Markov chains. In particular, we will consider to nearest-neighbor simple random walks.
%Both classes can be understood as possible stochastic generalizations of the one sided shift on $\mathbb{Z}_+$.

We say that $(X,T)$ is {\it hypercyclic}, or topologically transitive, if it has a dense orbit in X. This notion is equivalent that for all non empty open subsets $U$ and $V$ of $X$, there exists an integer $n \geq 0$ such that $T^{n}(U) \cap V$ is not empty. If moreover for every non-empty open set $V \subset  X$, the set $N(x, V) = \{k \in \mathbb{N},\; T^{k}(x) \in V \}$ has positive lower density, i.e
$
\liminf_{n \rightarrow \infty} \; \frac{1}{n} card (N(x, V) \cap [1, n]) > 0 \, ,
$
then we call $(X, T)$ {\it frequently hypercyclic}.
On the other hand, $(X, T)$ is said to be {\it supercyclic} if there exists $x \in X$ such that the projective orbit of $x$ is  is dense in the sphere $S^{1}= \{z \in X, \| z \|=1\}$, that is the set  $\{\lambda T^{n}(x),\; n \in \mathbb{N}, \lambda \in \mathbb{C}\}$ is dense  in X. We call  $(X, T)$  {\it Devaney chaotic} if it is hypercyclic, has a dense set of periodic points and has a sensitive dependence  on the initial conditions.

The study of those four properties  is a central problem in area of linear dynamical systems (see for instance \cite{bm} and \cite{em}).
Notice that the above properties can be studied in the context of more general topological space $X$ called Frechet spaces (the topology is induced by a sequence of semi-norms).

There are many examples of hypercyclic linear operators (see \cite{bm}) as the derivative operator on the Frechet space $H(\mathbb{C})$  of holomorphic maps on $\mathbb{C}$ endowed  with the topology of uniform convergence on compact sets, translation operator on  $H(\mathbb{C})$, classes of weighted shift operators acting on $X \in \{c_0, l^p, p \geq 1\}$.

However, the set of hypercyclic linear operators is  small. In fact it is proved that this set is nowhere dense in the set of continuous linear operators with respect to the norm topology (see \cite{bm}). An example of non hypercyclic operator is  the shift operator $S$ acting on $X \in \{c_0, l^q, q \geq 1\}$. This come from the fact that the norm of $S$ is less or equal to $1$. However, the shift operator is supercyclic and moreover for any $\lambda >1,\; \lambda S$ is frequently hypercyclic and chaotic (see \cite{bm}).

Here we are interested in operators associated to stochastic infinite matrices acting on a separable Banach space $X \in \{c_0, c, l^q, q \geq 1\}$.
In particular, we prove  that if  $A$ is a transition operator on an irreducible Markov chain with countable state space acting on $c$, then $A$ is not supercyclic. The result remains valid if we replace $c$ by $c_0$ or $l^q,\; q \geq 1$ in the positive recurrent case. A natural question is: what happens when the Markov chain is null recurrent or transient if $X \in \{c_0, l^q, q \geq 1\}?$ In order to study the last question, we consider
 transition operators $W_p$ (resp. $\overline{W}_p$) of nearest-neighbor simple asymmetric random walks on $\mathbb{Z}_+$ (resp. on $\mathbb{Z}$) with jump probability $p\in (0,1)$.

For the simple  asymmetric random walk on $\mathbb{Z}^{+}$ defined in $X \in \{c_0, c, l^q, q >1\}$, we prove that
if the random walk is transient ($p >1/2$), then $W_p$ is supercyclic and moreover  for all $\vert \lambda \vert> \frac{1}{2p -1},\; \lambda W_p$ is frequently hypercyclic and chaotic.
If  the random walk is null recurrent ($p=1/2$) and $X= l^1$, then $W_p$ is not supercyclic.

For the simple  asymmetric random walk on $\mathbb{Z}$, we prove that if $p \ne 1/2$ (transient case), $\lambda\overline{W}_p$ is not hypercyclic for all $\vert \lambda \vert > \frac{1}{\vert 1 -2p \vert}$.

We also consider transition operators spatially inhomogeneous simple random walks on $\mathbb{Z}_{+}$, that is operators $G_{\bar{p}} := G$ associated to a sequence of probabilities $\bar{p} = (p_{n})_{n \geq 0}$   and defined by $G_{0,0}= 1- p_0,\; G_{0,1}= p_0$ and for all $i \geq 1,\;  G_{i, j} =0 $ if $j \not \in \{ i-1, i+1 \}$, $ G_{i, i-1}= 1-p_i$ and $ G_{i, i+1}= p_i$.

In particular, we prove the following result:
Consider the sequence
$$
w_n= \frac{(1-p_1)(1-p_3)...(1-p_{n-1})}{p_1 p_3... p_{n-1}} \textrm{ for } n \textrm{ even},
$$
and
$$
w_n= \frac{(1-p_0)(1-p_2)...(1-p_{n-3})(1-p_{n-1})}{p_0 p_2... p_{n-3}p_{n-1}} \textrm{ for } n \textrm{ odd}.
$$

The following results hold:

1. If  $X= c_0$ and $\lim w_n=0$ or $X= l^q,\; q \geq 1$ and $\sum_{n=1}^{+\infty} w_n^q <+\infty$, then $G$  is supercyclic on $X$.

2. Let $X \in \{c_0, l^q,\; q \geq 1\}$ and   assume that there exists $\alpha >0$ such that $p_n \geq \frac{1}{2}+ \alpha$ for all $n \geq n_0$,   then there exists $\delta >1$ such that  $ \lambda G$ is frequently  hypercyclic and Devaney chaotic for all $\vert \lambda \vert >\delta$.

The last two  results can be extended for the spatially inhomogeneous simple random walks on $\mathbb{Z}$.

As a consequence of our dynamical study of random walks, we deduce that
if the Markov chain is null recurrent, it cannot be supercyclic on $l^1$ (see Proposition 4.6) or supercyclic on $c_0$ (see remark 4.4). We also deduce that, when the Markov chain is transient, it can have nice dynamical properties as supercyclicity, frequently hypercyclicity and  chaoticity on $X$ in $\{c_0, l^p, p \geq 1\}$ (see Theorems 4.1 and 4.8). We wonder if it is possible to construct transient Markov chains on $\mathbb{Z}_+$ or $\mathbb{Z}$ that are not supercyclic.

\smallskip

The paper is organized as follows: In section 2, we give some definitions and classical results. Section 3 describes the study of dynamical properties of Markov chain operators.
In section 4, we consider operators associated to the simple  asymmetric
  and also the spatially inhomogeneous simple random walks on $\mathbb{Z}_{+}$ and $\mathbb{Z}$.

%%%%%%%%%%%%%%%%%%%%%%%%%%%%%%%%%%%%%%%%%%%%%%%%%%%%%%%%%%%%%%%%%%%%%%%%%%%%%%%%%%%%%%%
%%%%%%%%%%%%%%%%%%%%%%%%%%%%%%%%%%%%%%%%%%%%%%%%%%%%%%%%%%%%%%%%%%%%%%%%%%%%%%%%%%%%%%%
%%%%%%%%%%%%%%%%%%%%%%%%%%%%%%%%%%%%%%%%%%%%%%%%%%%%%%%%%%%%%%%%%%%%%%%%%%%%%%%%%%%%%%%
%%%%%%%%%%%%%%%%%%%%%%%%%%%%%%%%%%%%%%%%%%%%%%%%%%%%%%%%%%%%%%%%%%%%%%%%%%%%%%%%%%%%%%%
%%%%%%%%%%%%%%%%%%%%%%%%%%%%%%%%%%%%%%%%%%%%%%%%%%%%%%%%%%%%%%%%%%%%%%%%%%%%%%%%%%%%%%%
%%%%%%%%%%%%%%%%%%%%%%%%%%%%%%%%%%%%%%%%%%%%%%%%%%%%%%%%%%%%%%%%%%%%%%%%%%%%%%%%%%%%%%%
%%%%%%%%%%%%%%%%%%%%%%%%%%%%%%%%%%%%%%%%%%%%%%%%%%%%%%%%%%%%%%%%%%%%%%%%%%%%%%%%%%%%%%%
%%%%%%%%%%%%%%%%%%%%%%%%%%%%%%%%%%%%%%%%%%%%%%%%%%%%%%%%%%%%%%%%%%%%%%%%%%%%%%%%%%%%%%%

\section{Definitions and classical results}

 To fix the notation we introduce here the proper definitions of the spaces mentioned above: Let $w= (w_n)_{n \geq 0}$ be a sequence of complex numbers. We put
$$
\| w \|_\oo =  \sup_{n \ge 0} |w_n| < \oo \, , \quad
\| w \|_q = \Big( \sum_{n\ge 0} |w_n|^q \Big)^\frac{1}{q} \, , \ 1 \le q < \oo \, ,
$$
and
$$
l^\oo = l^\oo(\mathbb{Z}_+) = \{ w \in \mathbb{C}^{\mathbb{Z}_+} : \| w \|_\oo < \oo \} \, ,
$$
$$
l^q = l^q(\mathbb{Z}_+) = \{ w \in \mathbb{C}^{\mathbb{Z}_+} : \| w \|_q < \oo \} \, ,
$$
$$
c = c(\mathbb{Z}_+) = \{ w \in l^\oo : w \textrm{ is convergent}  \} \, ,
$$
$$
c_0 = c_0 (\mathbb{Z}_+) = \{ w \in c : \lim_{n\rightarrow \oo} w_n = 0 \} \, .
$$

Now recall the definitions from the introduction. Other definitions related to linear dynamics will be needed.

\begin{definition}
Let $f : Y \to Y$ be a continuous map acting on some metric space
$(Y, d)$. We say that $f$ is  {\it Devaney chaotic} if

(1) $f$ is hypercyclic;

(2) $f$ has a dense set of periodic points;

(3) $f$ has a sensitive dependence on initial conditions: there exists $\delta > 0$ such that,
for any $x$ in $Y$ and every neighborhood $U$ of $x$, one can find $y  \in U$ and an
integer $ n >0$ such that $d(f^{n}(x), f^{n}(y)) \geq \delta$.
\end{definition}

\medskip

Let $X$ be a Banach separable space on $\mathbb{C}$ and $T: X \to X$ be a linear bounded operator on X.

\begin{definition} (see \cite{bm}).
We say that $T$ satisfies the hypercyclicity (resp. supercyclicity) criterion if there exists an increasing sequence of nonnegative integers $(n_k)_{k \geq 0}$, two dense subspaces of $X,\;  D_1$ and $D_2$ and a sequence of maps $S_{n_{k}}: D_2 \to X$ such that
\begin{enumerate}
\item
$\lim_{k} T^{n_{k}}(x)= \lim_{k} S_{n_{k}}(y)=0$
  (resp. $\lim_{k} \|\vert T^{n_{k}}(x)\| \| S_{n_{k}}(y)\| =0$),\; $ \forall x \in D_1,\; y \in D_2$.
\item
$\lim_{k} T^{n_{k}} \comp S_{n_{k}}(x)=x, \; \forall x \in D_2$.
\end{enumerate}
\end{definition}

\medskip

\begin{theorem}
The following properties are true
\begin{enumerate}
\item
If $T$  satisfies the hypercyclicity criterion, then $T$ is hypercyclic.
\item
$T$  satisfies the hypercyclicity criterion if and only if $T$ is topologically weakly mixing, i.e $T \times T$ is topologically mixing.
\item
If $T$  satisfies the supercyclicity criterion, then $T$ is supercyclic.
\end{enumerate}
\end{theorem}

\medskip

There is an efficient criterion that guarantees that $T$ is Devaney  chaotic and frequently hypercyclic (see \cite{bm}).

\medskip

\begin{theorem}
\label{crifre}
 Assume that there
exist a dense set $D \subset X$ and a map $S :D  \to D$ such that

\medskip

\begin{enumerate}
\item
For any $x \in D$, the series $\sum_{n=0}^{+\infty} T^{n}(x)$ and $\sum_{n=0}^{+\infty} S^{n}(x)$ are unconditionally convergent (all subseries of both series are convergent).
\item
For every $x \in D,\;  T \circ S(x)= x$,
\end{enumerate}
then $T$ is chaotic and frequently hypercyclic.
\end{theorem}

Concerning the dynamical properties of a linear dynamical system $(X,T)$, the spectrum of $T$ plays an important role.
We  denote by  $\sigma(X, T)$, $\sigma_{pt}(X,T)$, $\sigma_r(X,T)$ and $\sigma_c(X,T)$ respectively the spectrum, point spectrum, residual spectrum and continuous spectrum of $T$. Recall that $\lambda$ belongs to $ \sigma(X,T)$ (resp.  $\sigma_{pt}(X,T)$) if $(S- \lambda I)$ is not bijective (resp. not one to one).
If $(S- \lambda I)$ is one to one and not onto, then $\lambda \in \sigma_r(X,T)$
if $(S- \lambda I) (X)$ is not dense in $X$, otherwise, we say that $\lambda \in \sigma_c(X,T)$. Below we also use the notation $X'$ and $T'$ to indicate respectively the topological dual space and the dual operator associated to $(X,T)$.

\begin{lemma}  (\cite{bm})
\label{spectrhyper}
Let $X$ be a Banach separable space on $\mathbb{C}$ and $T: X \to X$ be a linear bounded operator on X.
\begin{enumerate}
\item
 If $T$ is hypercyclic then every connected component of the spectrum intersects the unit circle.
 \item
 If  $T$ is  hypercyclic, then $\sigma_{pt}(X',T')= \emptyset$,
 \item
If $T$ is supercyclic then there exists a real number $R \geq 0$ such every connected component of the spectrum intersects the  circle $\{z \in \mathbb{C},\; \vert z \vert = R\}$.
\item
If  $T$ is  supercyclic, then $\sigma_{pt}(X',T')$ contains at most one point.
\end{enumerate}
\end{lemma}

\medskip

In this paper, we will use only items 2) and 4) of Lemma \ref{spectrhyper}. 1) and 3) are used in \cite{acmv} for the study of dynamical properties of Markov chains associated to stochastic adding machines.
\medskip

\begin{remark}
If $T$ is not supercyclic then  $\lambda T$ is not hypercyclic for every fixed $\lambda$. However, it is possible to have $T$ supercyclic and $\lambda T$ not hypercyclic for all sufficiently large (but fixed) $\lambda$.
\end{remark}

\medskip

%%%%%%%%%%%%%%%%%%%%%%%%%%%%%%%%%%%%%%%%%%%%%%%%%%%%%%%%%%%%%%%%%%%%%%%%%%%%%%%%%%%%%%%
%%%%%%%%%%%%%%%%%%%%%%%%%%%%%%%%%%%%%%%%%%%%%%%%%%%%%%%%%%%%%%%%%%%%%%%%%%%%%%%%%%%%%%%
%%%%%%%%%%%%%%%%%%%%%%%%%%%%%%%%%%%%%%%%%%%%%%%%%%%%%%%%%%%%%%%%%%%%%%%%%%%%%%%%%%%%%%%
%%%%%%%%%%%%%%%%%%%%%%%%%%%%%%%%%%%%%%%%%%%%%%%%%%%%%%%%%%%%%%%%%%%%%%%%%%%%%%%%%%%%%%%

\section {Dynamical properties of Markov Chains Operators}

\medskip

Let $Y = (Y_n)_{n\ge 1}$ be a discrete time irreducible Markov chain with countable state space $E$ and with transition operator $A = [A_{i,j}]_{i,j\in E}$ (irreducible means that for each pair $i$, $j \in E$ there exists a nonnegative integer  $n$ such that $A^n_{i.j} > 0$).
%on a given probability space $(\Omega,\mathcal{F},P)$.
The Markov chain $Y$ is said to be recurrent if the probability of visiting any given state is equal to one, otherwise $Y$ is said to be transient. The Markov chain $Y$ is called positive recurrent if it has an invariant probability distribution, i.e., there exists $u \in l^1$ such that $uA = u$. Every positive recurrent Markov chain is recurrent. If $Y$ is recurrent but not positive recurrent, it is called null recurrent.

For the transient and null recurrent cases we have the following well-known equivalent definitions (see \cite{ross}):
\begin{enumerate}
\item[(i)] $A$ is transient if and only if $\sum_{n= 1}^{+\infty} A^n_{i,j} < \infty$ for all $i,j \in E$.
\item[(ii)] $A$ is null recurrent if and only if $\lim_{n \rightarrow \infty} A^n_{i,j} = 0$ and $\sum_{n= 1}^{+\infty} A^n_{i,j} = \infty$ for all $i,j \in E$.
\end{enumerate}

\bigskip

\begin{proposition}
\label{c-ciclic}
Let $A$ be a transition operator on an irreducible Markov chain with countable state space acting on $c$, then $A$ is not supercyclic. The result remains valid if we replace $c$ by $c_0$ in the positive recurrent case.
%Let $A$ be a transition operator  Markov Chain with countable state space acting on $c$ \textcolor{red}{or $c_0$}, then $\lambda A$ is not hypercyclic for all $ \lambda \in \mathbb{C}$.
%Moreover if the Markov chain is transient or null recurrent, then $A$ is not supercyclic on $c$.
\end{proposition}

\noindent {\bf Proof:}
Consider an enumeration of the state space so that we can consider $E = \mathbb{N}$ and the stochastic matrix $A = [A_{ij}]_{i,j \in \mathbb{N}}$ associated with the transition operator $A$. Assume that $A$ is is transient or null recurrent, then $\lim_{n\rightarrow \infty} A^{n}_{i,j} = 0$ for every $i$ and $j$.

Now fix $y \in c-c_0$ (we do not need to consider the case $y \in c_0$ while considering density of orbits of $y$ under $A$ or $\lambda A$ because $c_0$ is a closed invariant subspace). Suppose that $\lim_{n \rightarrow \infty} y_n = \alpha \in \mathbb{C} - \{0\}$. We have that
\begin{equation}
\label{eq:convTn1}
\lim_{n\rightarrow \infty} \big( A^{n}y \big)_i  = \alpha \, ,
\end{equation}
for every $i \in \mathbb{N}$. Indeed, since $\sum_{j=1}^{+\infty} A^{n}_{i,j} = 1$, for every $n \in \mathbb{N}$
$$
\Big| \big( A^{n}y \big)_i - \alpha \Big| = \Big| \sum_{j=1}^{+\infty} A^{n}_{i,j} (y_j - \alpha) \Big| \le
\big( \|y\|_{\infty} + \vert \alpha \vert \big) \sum_{j=1}^m A^{n}_{i,j} + \sup_{j \ge m+1} |y_j - \alpha| \, .
$$
The second term in the rightmost side of the previous expression can be made arbitrarily small by choosing $m$ sufficiently large while the first one goes to zero as $n$ tends to $+\infty$ for every choice of $m$. Hence \eqref{eq:convTn1} holds.

From \eqref{eq:convTn1}, we have that
$$
\liminf_{n \rightarrow \infty} \frac{|(A^n y)_i|}{\| A^n y \|_{\infty}} \ge \frac{\vert \alpha \vert}{\| y \|_{\infty}} > 0 \, ,
$$
and then $\{ (A^n y)/\| A^n y \|_{\infty} : n \ge 1 \}$ is not dense in the unit sphere of $c$ centered at $0$ which implies that $\{ \lambda A^n y \, : \ \lambda \in \mathbb{C}, \ n \ge 1 \}$ is not a dense subset of $c$. Since $y$ is arbitrary, $A$ is not supercyclic on $c$.

Now, assume that $A$ is positive recurrent, then there exists an invariant measure $u \in l^1 \setminus \{0\}$ such that $u A= u$, hence $u A^n= u$ for all integer $n \geq 1$. Suppose that $A$ is supercyclic. Take $y \in c \cap S^1_c$, where $S^1_c = \{ x \in c : \| x \|_\infty = 1  \}$ such that projective orbit of $y$ under $A$ is dense in $S^1_c$, then for all $x \in S^1_c$, there exists an increasing  sequence $(n_k)_{k \geq 0}$ such that $\lim_{k \rightarrow \infty} \frac{A^{n_k}y}{\| A^{n_k} y \|_{\infty}} = x$. Since $\| A^{n_k} y \|_{\infty} \le \| y \|_{\infty} = 1$, then
$$ |<u, x>| = \lim_{k} \frac{ | < u , A^{n_k} y > |}{\| A^{n_k} y \|_{\infty} } = \lim_{k} \frac{ |< u , y >| }{\| A^{n_k} y \|_{\infty} } \ge
| < u , y > |.$$
Where $<u, z>$ is the scalar product between $u$ and $z$ for $z$ in $c$.
Since $x$ is arbitrary in $S^1_c$, the last inequality implies that $u=0$ and this is an absurd. Hence the projective orbit of $y$ under $A$ could not be dense in $S^1_c$ which means that $A$ is not supercyclic. \\
To finish the proof we just point out that in the positive recurrent case, the same proof holds if we replace $c$ by $c_0$. $\square$

%Suppose that there exists $\lambda \in \mathbb{C}$ such that $\lambda A$ is hypercyclic. Take $x \in c$ such the orbit of $x$ under $A$ is dense in $X$, then for all $y \in X$, there exists an increasing  sequence $(n_k)_{k \geq 0}$ such that $\lim A^{n_k}x = y$. Hence $$ <u, y> = \lim _{k}  u A^{n_k}x  = <u, x>.$$Hence $u=0$ and this is an absurd.

\medskip

%{\bf Question:} Does there exist a supercyclic transition operator $A$ on $c$ associated to a positive recurrent Markov chain?

%\begin{proposition} \label{posit}
%Let  $A: X \to X$ where $X \in \{l^q,\; q \geq 1\}$ be the transition operator of a positive recurrent stochastic Markov chain,  then $\lambda A$ is not hypercyclic for all $ \lambda \in \mathbb{C}$.\end{proposition}

%\noindent {\bf Proof:}
%Since $A$ is positive recurrent, then there exists an invariant measure $u \in l^1$ such that $u A= u$, hence for all $\lambda \in \mathbb{C},\; \lambda \in \sigma_{pt} ((\lambda A)', l^1)$. Since the dual of  $c$ is $l^1$ and the dual of $l^q,\; q \geq 1$ is $l^{\frac{q}{q-1}}$ which contain in $l^1$, we deduce that $\lambda$ belongs to $\sigma_{pt} ((\lambda A)', X')$ where $X= c$ or $X= l^q.$ Using item (2) of Lemma \ref{spectrhyper}, we obtain the result. \hfill $\Box$

\medskip

Another result is:

\begin{proposition}
\label{lq-hyp-c0}
 Let $q \geq 1$ and $A: l^q \to l^q$ be an hypercyclic (supercyclic) operator on $l^q$, then
 \begin{enumerate}
 \item
  If  $A(c_0) \subset c_0$,  then $A$ is also hypercyclic (supercyclic) on $c_0$.
  \item
  If $r > q$ and $A(l^r) \subset l^r$, then $A$ is also hypercyclic (supercyclic) on $l^r$.
 \end {enumerate}
\end{proposition}

\noindent {\bf Proof:}
(1) Suppose that $A$ is hypercyclic   on $l^q$ and let $x \in l^q$ be a hypercyclic vector,\; i.e $\overline{ O(x)}= l^q$. Now fix $y \in c_0$ and $\epsilon >0$. Take $m \in \mathbb{N}$ such that $\sup_{i > m} |y_i| \le \epsilon/2$. Define $y^{(m)}$ as
$$
y^{(m)}_i = \left\{
\begin{array}{cl}
y_i &, \ 1\le i \le m ,\\
0 &, \ \textrm{otherwise},
\end{array}
\right.
$$
for every $i \in \mathbb{N}$. Since $y^{(m)} \in l^q$, there exists $n \in \mathbb{N}$ such that $\| A^nx - y^{(m)} \|_\infty \le \| A^nx - y^{(m)} \|_q \le \epsilon/2$. Therefore
$$
\| A^nx - y \|_\infty \le \| A^nx - y^{(m)} \|_\infty + \| y^{(m)} - y \|_\infty \le \epsilon \, .
$$
Since $\epsilon$ and $y$ are arbitrary, $x$ is a hypercyclic vector in $c_0$.

The proof in the supercyclic case is analogous.

(2) The  proof is analogous to item 1) and come from the fact that if $1 \leq q <r$, then $l^q \subset l^r$.
$\square$

\medskip

\begin{corollary}
\label{posit}
Let  $A: X \to X$ where $X \in \{l^q,\; q \geq 1\}$ be the transition operator of a irreducible positive recurrent stochastic Markov chain, then $A$ is not supercyclic.
\end{corollary}

\noindent {\bf Proof:}
From Proposition \ref{c-ciclic} we have that $A$ acting on $c_0$ is not supercyclic. Thus from (1) in Proposition \ref{lq-hyp-c0} we obtain that $A$ acting on $X$ is not supercyclic.
$\square$

\medskip

\begin{remark}
Since an operator $A$ is supercyclic if and only if $cA$ is supercyclic for $c \neq 0$, then all the previous results in this section hold for operators associated to countable non-negative irreducible matrices with each line having the same sum of their entries.
\end{remark}

\medskip

\noindent \textbf{Question:}
What happens if $A:X \to X$ is a countable infinite non-negative irreducible  matrix where the the sum of entries of lines is not constant?

Is $A$ not supercyclic on $c$?

If $A$ is positive recurrent (see \cite{k}  for the definition), Can we prove that $A$ is not supercyclic in $X \in \{c_0, c,  l^q,\; q \geq 1\}$?

\section{Simple Random Walks}

\medskip

Consider the nearest neighbor simple random walk on $\mathbb{Z}_+$ with partial reflection at the boundary and jump probability $p\in (0,1)$ (when at zero, the walk stays at zero with probability $1-p$). Denote by $W_p := W= (W_{i,j})_{i, j \geq 0}$ its transition operator. We have $ W_{0,0}= 1-p,\; W_{0,1}= p$ and for all $i \geq 1,\;  W_{i, j} =0 $ if $j \not \in \{ i-1, i+1 \},\; W_{i,i-1}= 1-p,\; W_{i,i+1}= p$ for all $i \geq 1$. We have

$$ W_p=
\tiny{
\left[
\begin{array}{cccccccccc}
1-p & p & 0   & 0   & 0    & 0   &  0  &  0  & 0  & \cdots  \\
1-p  & 0 & p  & 0 & 0  &  0  & 0  & 0  & 0  & \cdots  \\
0  & 1-p  & 0  & p & 0 & 0 & 0 & 0 & 0 & \cdots  \\
0 & 0  & 1-p & 0 & p & 0 & 0 & 0 & 0 & \cdots  \\
\vdots & \vdots & \vdots & \vdots & \vdots & \vdots & \vdots & \vdots & \vdots  & \ddots
\end{array}
\right]}
$$

It is known (\cite{ross}) that the simple random walk on $\mathbb{Z}_{+}$ is positive recurrent if $p<1/2$, null recurrent if $p=1/2$ and transient if $ p>1/2$.

In particular, from Proposition \ref{c-ciclic} and Corollary \ref{posit}, we have that $W_p$ acting on $X \in \{c_0,c,l^q \ q\ge 1\}$ is not supercyclic if $p<1/2$.

\medskip

\begin{theorem}
\label{passeio}
Let $X \in \{c_0, l^q,\; q \geq 1\}$. If $p >1/2$, then the infinite matrix $W_p$ of the simple random walk on $\mathbb{Z}_{+}$ is supercyclic on $X$. Moreover $\lambda W_p$ is frequently  hypercyclic and chaotic for all $\vert \lambda \vert >\frac{1}{2p-1}$.
\end{theorem}

\smallskip

Before we prove Theorem \ref{passeio} we need three technical results:

\begin{lemma}
\label{zer}
Let $X \in \{c_0, l^q,\; q \geq 1\}$, then $\sigma_{pt}(X,W_p)$ is not empty if and only if $p >1/2$, moreover, in this case $0 \in \sigma_{pt}(X,W_p)$.
\end{lemma}

\noindent {\bf Proof:}
Let $ \lambda $ be an element of $\sigma_{pt}(W_p)$ and $u= (u_n)_{n \geq 0}$ be an eigenvector associated to $\lambda$, then
$$ (1- p- \lambda) u_0+ p u_1=0,\; (1-p)u_n - \lambda u_{n+1}+ p u_{n+2}=0,\; \forall n \geq 0.$$
We deduce that there exists a sequence of complex numbers $(q_{n})_{n \geq 0}$  where $q_0=1,\; q_1= \frac{\lambda + p-1}{p}$ and  $u_n= q_n u_0$.
Moreover
$$\begin{bmatrix}q_n\\ q_{n-1}\end{bmatrix} = M \begin{bmatrix}q_{n-1}\\ q_{n-2}\end{bmatrix},\;\; \forall n \geq 2.$$
Where $M= \begin{bmatrix}\frac{\lambda} { p}& \frac{p-1}{p}&\\ 1&0\end{bmatrix}$.
Hence
$\begin{bmatrix}q_n\\ q_{n-1}\end{bmatrix} = M^{n-1} \begin{bmatrix}q_1\\ q_{0}\end{bmatrix}$ for all $n \geq 2$.
Assume that $\lambda^2 \neq 4 p (1-p)$, then the matrix $M$ have distinct eigenvalues and hence it is diagonalizable. Therefore, there exist $c, d \in \mathbb{C}\setminus \{0\}$ such that $q_n= c \alpha^n + d \beta^n$ for all integer $n \geq 0$, where $\alpha, \beta$ are the eigenvalues of $M$.
Since $\alpha \beta= det (M)= \frac{1-p}{p}$, then if $0 <p <1/2$, we have  $ \alpha \beta > 1$ . Then either $\vert \alpha \vert > 1$ or  $\vert \beta \vert > 1$,
Hence $(q_n)_{n \geq 0}$ is not bounded and therefore  the point spectrum of $W_p$ is  empty.
If $p=1/2$, then either ($\vert \alpha \vert > 1$ or  $\vert \beta \vert > 1$) or $  \alpha, \beta$ are conjugated complex numbers of modulus $= 1$. In both cases the point spectrum of $W_p$ is  empty.

If  $\lambda^2 = 4 p (1-p)$, then the matrix $M$  is not diagonalisable and has a unique eigenvalue $\theta$. In this case, there exist $e,  f \in \mathbb{C} \setminus \{0\}$ such that $q_n= (e+ f n) \theta^n$ for all integer $n \geq 0$.
If $p  \leq 1/2$, then $\vert \theta \vert = \sqrt {\frac{1-p}{p} } \geq 1$, hence $q_n$ is not bounded.
 We deduce that  the point spectrum of $W_p$ is  empty.

Now assume that $1/2 <p <1$ and $\lambda= 0$ ($M$ diagonalizable), then $\alpha, \beta \in  \{- \sqrt { \frac{1-p}{p}} i,\; \sqrt { \frac{1-p}{p}} i\}$, therefore $\alpha, \beta$ are complex conjugated numbers of modulus $< 1$.
Hence $0 \in \sigma_{pt}(X,W_p)$.
If $p=1$ and $\lambda= 0$, then $q_0=e$ and $q_n=0$ for all $n \geq 1$. Thus $0 \in \sigma_{pt}(X,W_p)$.
$\square$

\medskip

\begin{remark}
\label{zeropoint}
Consider $M$ as in the proof of Lemma \ref{zer}.
1. Since the eigenvalues of $M$ depend continuously of $\lambda$, we deduce that for all $p > 1/2$, there exists $0 < r_p < 1$ such that $D(0, r_p) \subset \sigma_{pt}(X,W_p)$. In particular, we can prove that $[0, 2 \sqrt{1-p})
\subset \sigma_{pt}(X,W_p)$.

2. If $p=1/2$  the eigenvalues of $M$ are $\lambda \pm \sqrt{\lambda^2 -1}$, we deduce that if $X= l^{\infty}$, the interval $ ]-1, 1[ \subset \sigma_{pt}(l^{\infty},W_p)$.
\end{remark}

\medskip

\begin{lemma}
\label{inversa1}
Let $X \in \{c_0, l^q,\; q \geq 1\}$ and  $v= (v_{i})_{i \geq 0} \in X$.
Let $a= (a_n)_{n \geq 0} \in l^{1}$ and $x= (x_n)_{n \geq 0}$ defined by
$$x_n= \sum_{k=0}^{n} a_k v_{n-k},\; \forall n \in \mathbb{Z}_{+}.$$ Then $(x_n)_{n \geq 0} \in X$, moreover
 $$\| x \| \leq \||a \|_1 \| v \|.$$
\end{lemma}

\noindent {\bf Proof:}
By putting $v_k= 0$ for all $k <-1$, we can assume that
$x_n=  \sum_{k=0}^{+\infty} a_k v_{n-k}$ for all $n \geq 0$.

Now suppose that $X= c_0$. For each $i \in \N$
$$
|x_n| \leq  \Bigg(\sum_{k=0}^i \vert a_k \vert \Bigg)
               \sup_{0 \leq k \leq i} \|v_{n-k}\|
           + \Bigg(\sum_{k=i+1}^\infty \vert a_k \vert \Bigg)
               \|v\|_{\infty}.
$$
Since $ (v_{n})_{n \geq 0} \in c_0$.
and $(a_n)_{n \geq 0} \in l^1$, we deduce that
$(x_n)_{n \geq 0} \in c_0$.

If $X= l^1$, we have
$$
\sum_{n =0}^{+\infty} |x_n|
  \leq  \sum_{n =0}^{+\infty}\sum_{k=0}^\infty \vert a_k v_{n-k} \vert
     \leq \Bigg(\sum_{k=0}^\infty \vert a_k \vert\Bigg)
         \Bigg(\sum_{n =0}^{+\infty} \vert v_{n} \vert \Bigg) =\|a \|_1  \|v \|_1 .
$$
Now, assume that $X= l^q,\; 1 < q < \infty$, we consider its conjugate exponent
$r$ (i.e., $1/q + 1/r= 1$).
We have $|x_n|
  \sum_{k=0}^\infty\vert a_k \vert^\frac{1}{r}\vert a_k \vert^{ (1- 1/r)}\vert v_{n-k} \vert$. Hence
 By  H\"older's inequality, we obtain
$$
|x_n|
  \leq   \Bigg(\sum_{k=0}^\infty\vert a_k \vert\Bigg)^\frac{1}{r}
       \Bigg(\sum_{k=0}^\infty \vert a_k \vert^{q (1- 1/r)}\vert v_{n-k} \vert^q
             \Bigg)^\frac{1}{q}.
$$
As a consequence,
$
\sum_{n =0}^{+\infty} \vert x_n \vert^q
  \leq \Bigg(\sum_{k=0}^\infty\vert a_k \vert\Bigg)^{q-1}
          \sum_{n =0}^{+\infty}\sum_{k=0}^\infty \vert a_k \vert \vert v_{n-k} \vert^q.$
        Hence
        $$
\sum_{n =0}^{+\infty} \vert x_n \vert^q
  \leq \Bigg(\sum_{k=0}^\infty\vert a_k \vert \Bigg)^{q}
        \sum_{n =0}^{+\infty} \vert v_{n} \vert^q.$$
$\square$

\begin{proposition}
\label{inversa}
Let $X \in \{c_0, l^q,\; q \geq 1\}$ and $p >1/2$
then for $v= (v_{i})_{i \geq 0} \in X$ there exists $u= (u_{i})_{i \geq 0} \in X$ such that $W(u) = v$. Moreover, $u$ is explicitly given by
\begin{eqnarray}
\label{frma}
u_{n}= \frac{1}{p} \; \Big( \sum_{j=0}^{\lfloor n/2 \rfloor} \Big( \frac{p-1}{p} \Big)^{j} v_{n-2j-1} + \Big( \frac{p-1}{p} \Big)^{\lfloor n/2 \rfloor + 1} p u_0 \Big)
\end{eqnarray}
where $\lfloor n/2 \rfloor $  is the largest integer  $< n/2$.
\end {proposition}

\noindent {\bf Proof:}
Fix $v= (v_{i})_{i \geq 0} \in X$ and $u= (u_{i})_{i \geq 0} \in l^{\infty}$ such that $W(u)= v$, then

\begin{eqnarray}
\label{fre} u_1= \frac{v_{0}}{p}+ \frac{p-1}{p} u_{0},\; u_n = \frac{v_{n-1}}{p}+ \frac{p-1}{p} u_{n-2},\; \forall n \geq 2.
\end{eqnarray}

Then we obtain by induction that for all integer $n \geq 2$,
$$
u_{n}= 1/p \; ( v_{n-1}+ \gamma v_{n-3} + \gamma^2 v_{n-5}+ \ldots \gamma^k v_{n-2k-1}+ \ldots \gamma^{\lfloor n/2 \rfloor} v_t + \gamma^{\lfloor n/2 \rfloor+1} p u_0 ),
$$
where $\gamma= \frac{p-1}{p}$, $t= 0$ if $n $ is odd and $t=1$ otherwise. Thus we have (\ref{frma}).

By Lemma \ref{inversa1}, we deduce that $u$ belongs to $X$.

Observe that if $u \in c_0$, then
\begin{eqnarray}
\label{fie}\| u \|_{\infty} \leq \delta \; max ( \| v \|_{\infty}, \vert u_0 \vert)  \mbox { where } \delta = \frac{1}{p}\sum_{n=0}^{+\infty} \gamma ^n= \frac{1}{2p-1}.
\end{eqnarray}

If $u \in l^q,\; q \geq 1$, then by Lemma \ref {inversa1},
 \begin{eqnarray}
\label{fie2} \| u \|_q \leq \frac{1}{2p-1} \| v \|_q + \vert u_0 \vert   \Big( \sum_{n=1}^{+\infty} \gamma^{nq} \Big)^{1/q}.
\end{eqnarray}

\begin{remark}
\label{zer}
If $u_0= 0$, then $\| u \|_q \leq \frac{1}{2p-1} \| v \|$ (in $X$).
\end{remark}

\bigskip

\begin{lemma}
\label{kernel}
Let $X \in \{c_0, l^q,\; q \geq 1\},\; p >1/2$, then
$D= \bigcup_{n=1}^{+\infty} Ker (W^n)$ is dense in $X$.
\end {lemma}

\noindent {\bf Proof:}
Fix $X$ in $\{c_0, l^q,\; q \geq 1\}$. By Proposition \ref{zer}, we have that $0 \in \sigma_{pt}(X,W)$. Then for all integer $n \geq 1,\; Ker (W^n) $ is not empty.

\smallskip

{\bf Claim:} {\em for all integer $n \geq 1$, there exists $V_{0,n},\ldots , V_{n-1,n} \in l^{\infty}$, linearly independent such that if $u=(u_{i})_{i \geq 0} \in  Ker (W^{n})$, then $ u = \sum_{i=0}^{n-1} u_i V_{i,n}$.}

\smallskip

Indeed, since $W_{j,k}=0$ for all $k \geq j+2$, we deduce that
for all integer $n \geq 2,\; W^{n}_{j,k}=0$  for all integer $k \geq j+n+1$.

Assume that $u=(u_i)_{i \geq 0} \in  Ker (W^n)$. The relation $\sum_{k=0}^{j+n} W^{n}_{j,k} \, u_k=0$ holds for all $j \in \mathbb{N}$. Since $W_{0,n}^{n} > 0$, we deduce that
$$u_{n}=  \sum_{i=0}^{n-1} u_i c_{i,n,n} \, ,$$
where $c_{i,n,n}= -\frac{W^{n}_{0,i}}{W^{n}_{0,n}}$
for all $i=0,\ldots, n-1$.
We also obtain by induction that
$u_{k}=  \sum_{i=0}^{n-1} u_i c_{i,k,n}$ for all $k \geq n,$ where $c_{i,k,n}$ are real numbers.

For all $i \in \{0,1,\ldots n-1\}$, define the infinite vector  $V_{i,n}= (V_{i,n}(k))_{k \geq 0}$ by putting $V_{i,n}(k)= c_{i,k,n}$ for all $k \geq n$  and $V_{i,n}(k)= \delta_{i,k}$ for all $0 \leq k <n$. Then, we obtain the claim.

\smallskip

Now observe that  $V_{i, n} \in X$ for every integer $n$ and $i=0,\ldots, n$. Indeed  for all $i=0,\ldots n$, we have $W^{n-1} V_{i, n} \in ker W$ that is contained in $X$. Hence by Lemma \ref{inversa}, we deduce that $W^{n-2} V_{i, n} \in X$ and   continuing by the same way, we obtain that  $V_{i, n} \in X$.

Now, let $z= (z_i)_{i \geq 0} \in X$, such that $z_i=0$ for all $i >n$ where $n$ is a large integer
number, then $z$ can be approximated by the vector $\sum _{i=0}^{n-1} z_i V_{i,n}$ which belongs to
the set $D$. Hence the $ D$ is dense in $X$. $\square$

\bigskip

\bigskip

{\bf Proof of Theorem \ref{passeio}:}
First we prove that $W$ is supercyclic. Recall the definition of $D$ from the statement of Lemma \ref{kernel}. By Proposition \ref{inversa}, for every $v \in D$, we can choose $Sv \in D$ such that $W(Sv)= v$. Using the fact that $Sv \in D$, we prove by induction that $W^n (S^n v)= v$ for all $v \in D$. On the other hand, since for all $u \in D$, there exists a nonnegative integer $N$ such that $W^{n}(u)=0$ for all $n \geq N$, we deduce that
 $\lim \| W^{n}(u) \| \, \| S^{n}(v)\| =0, \; \forall u, v \in D$. Hence $W$ satisfies the supercyclicity criterion. Thus $W$ is supercyclic.

\vspace{1em}

\noindent {\bf Claim:  $\lambda W$ is frequently  hypercyclic and chaotic for all $\vert \lambda \vert> \frac{1}{2p-1} $.}

Indeed, let $ v= (v_{i})_{i \geq 0} \in X$ and $u= (u_{i})_{i \geq 0} \in l^{\infty}$ such that $W(u)= v$, then $u$ satisfies (\ref{frma}).
 Putting $u_0=0$, we obtain $S(v)= (0, u_1,u_2 \ldots)$ and $W(Sv)= v$.

We also have  by remark \ref{zer} that $\| Sv \| \leq \frac{1}{2p-1} \| v \|$.

On the other hand, since $S(v)_0= 0$, we obtain by (\ref{fre}) that
$$S^2(v)= (0, 0, (S^2 v)_2, (S^2 v)_{3},\ldots).$$
We deduce that  for all integer $n \geq 0 $
$$ S^n(v)= (\underbrace{0,\ldots, 0}_n, (S^n v)_n, (S^n v)_{n+1},\ldots) \mbox { and } \| S^n (v) \| \leq \Bigg(\frac{1}{2p-1}\Bigg)^{n} \| v \|.$$
Let $\lambda$ be a complex number such that $\vert \lambda \vert > \delta$, then $\| \lambda^{-n} S^n (v) \|$ converges to $0$ exponentially as $n$ goes to $+\infty$.

Taking $W'= \lambda W$ and $S' = \lambda^{-1} S$ and $D= \cup_{n=0}^{+\infty} Ker(W^n)$, we obtain that the series $\sum_{n=0}^{+\infty} W'^{n}(x)$ and $\sum_{n=0}^{+\infty} S'^{n}(x)$ are absolutely convergent and hence unconditionally convergent for all $x \in D$, moreover
 $ W' \circ S'= I$ on $D$, then we are done by Theorem \ref{crifre}. $\square$

\medskip

%\begin{corollary}If  $p <1/2$, then $\lambda W_p$ is not hypercyclic in $l^1$ for all $\vert \lambda \vert >1$.
%\end{corollary}

%\medskip

%Now, we will prove a stronger result:

%\begin{theorem}Let $X \in \{c_0, c, l^q,\; q > 1\}$ and   $p <1/2$, then the transition operator $W_p$  associated to the asymmetric random walk on $\mathbb{Z}_{+}$ is not supercyclic.\end{theorem}

%\noindent{\bf{Proof:}}
%We have for all $i \geq 1,\; W^{*}_{i, i-1}= p,\; W^{*}_{i, i+1}= p$ and $T^{*}_{0,0}= W^{*}_{0,1}= 1-p$, we deduce that if $ \lambda $ is an element of $\sigma_{pt}(W^{*})$ and $u= (u_n)_{n \geq 0}$ an eigenvector associated to $\lambda$,
%then for all $n \geq 2$ such that $u_n= l_n u_0$, where $l_0=1,\; l_1= \frac{\lambda + p-1}{1-p}$. and
%$\begin{bmatrix}l_n\\ l_{n-1}\end{bmatrix} = N^{n-1} \begin{bmatrix}q_1\\ q_{0}\end{bmatrix},\;\; \forall n \geq 2,$ where $N= \begin{bmatrix}\frac{\lambda} { 1-p}& \frac{p}{p-1}&\\ 1&0\end{bmatrix}$. We deduce that for all $ p < 1/2$, there there exists $0 < r_p < 1$ such that $D(0, r_p) \subset \sigma_{pt}(W_p^{*})$. Using Lemma (\ref{pointsp}), we obtain the result. $\square$

\medskip

\begin{proposition} \label{notSRW}
For $p =1/2$, the operator $W_p$ acting on $l^1$ is not supercyclic
\end{proposition}

\noindent \textbf{Proof:} Note that $W_p$ is symmetric for $p=1/2$, then by (2) in remark \ref{zeropoint} we have that $\sigma_{pt} ((l^1)',W_p') = \sigma_{pt} (l^\infty,W_p) \supset ]-1,1[$. By (4) in Lemma (\ref{spectrhyper}), we obtain the result. $\square$

\bigskip

\noindent {\bf Question:}
For $p=1/2$, is the operator $W_p$ acting on $c_0$ or $l^q,\; q >1$ not supercyclic?

\subsection {Simple Random Walks  on $\mathbb{Z}$}

Consider the simple random walk on $\mathbb{Z}$ with jump probability $p\in (0,1)$, i.e, at each time the random walk jumps one unit to the right with probability p, otherwise it jumps one unit to the left.
Denote by $\overline{W}_p:=\overline{W}$ its transition operator. For all $i, j \in \mathbb{Z}$,  We have $\overline{W}_{i,j}=0$ if $j \ne i-1$ or $j \ne i+1$ and $\overline{W}_{i, i-1}=1-p,\; \overline{W}_{i, i+1}=p$.

The simple random walk on $\mathbb{Z}$ is null recurrent if $p=1/2$, otherwise it is
transient.

\begin{proposition}
 If $p \ne 1/2$, then $\lambda \overline{W}_p$ is not hypercyclic on $X \in \{c_0, l^q,\; q \geq 1\}$, for all $\vert \lambda \vert \geq \frac{1}{\vert 1- 2p  \vert}$.
  If $p =1/2$, then  $\overline{W}_{p}$  is not supercyclic on $l^1$.
\end{proposition}

\noindent \textbf{Proof:}
Let $X \in \{c_0, l^q,\; q \geq 1\}$ and $x= (x_{i})_{i \in \mathbb{Z}} \in X$, then $S(x) = (1-p) y+ p z$ where $y= (y_i)_{i \in \mathbb{Z}}$ and $z= (z_i)_{i \in \mathbb{Z}}$ satisfy $y_i= x_{i-1}$ and $z_i= x_{i+1}$ for all $i$.

Hence
$$\| S(x) \| \geq  (1-p) \| y \| - p \| z \|.$$
Since $\| y \| = \| z \|= \| x \|$,
we deduce that $\| S(x) \| \geq \vert 1 -2p \vert \|  x \|$. Hence
 $$\| S^n(x) \| \geq \vert 1 -2p \vert^n \|  x \| \mbox {  for all } n\ge 1.$$
Then $\lambda \overline{W}_p$ is not hypercyclic on $X \in \{c_0, l^q,\; q \geq 1\}$, for all $\vert \lambda \vert \geq \frac{1}{\vert 1- 2p  \vert}$.

Now, assume that $p=1/2$. Note that $\overline{W}_p$ is symmetric, then by (2) in remark \ref{zeropoint} we have that $\sigma_{pt} ((l^1)',W_p') = \sigma_{pt} (l^\infty,W_p)$.
We can prove that $]-1, 1[ \subset \sigma_{pt} (l^\infty,W_p)$.  By (4) in Lemma (\ref{spectrhyper}), we obtain  $\overline{W}$  is not supercyclic on $l^1$.
$\square$

\medskip

\begin{remark}
Dynamical properties of simple Random Walks  on $\mathbb{Z}_{+}$ and
on $\mathbb{Z}$ are different in case where they are transient.
\end{remark}

\textbf{Question:} Is $\lambda \overline{W}_p$ not hypercyclic on $X \in \{c_0, l^q,\; q \geq 1\}$ for all $|\lambda| \ge 1$? Can $\overline{W}_p$ be supercyclic?

\

%%%%%%%%%%%%%%%%%%%%%%%%%%%%%%%%%%%%%%%%%%%%%%%%%%%%%%%%%%%%%%%%%%%%%%%%%%%%%%%%%%
%%%%%%%%%%%%%%%%%%%%%%%%%%%%%%%%%%%%%%%%%%%%%%%%%%%%%%%%%%%%%%%%%%%%%%%%%%%%%%%%%%
%%%%%%%%%%%%%%%%%%%%%%%%%%%%%%%%%%%%%%%%%%%%%%%%%%%%%%%%%%%%%%%%%%%%%%%%%%%%%%%%%%
%%%%%%%%%%%%%%%%%%%%%%%%%%%%%%%%%%%%%%%%%%%%%%%%%%%%%%%%%%%%%%%%%%%%%%%%%%%%%%%%%%

\subsection {Spatially inhomogeneous simple random walks on $\mathbb{Z}_{+}$}

In this section we consider spatially inhomogeneous simple random walks on $\mathbb{Z}_{+}$, or discrete birth and death processes. Let $\bar{p} = (p_{n})_{n \geq 0}$ be a sequence of probabilities, the simple random walk on $\mathbb{Z}_{+}$ associated to $\bar{p}$ is a Markov chain with transition probability $G_{\bar{p}} := G$ defined by $G_{0,0}= 1- p_0,\; G_{0,1}= p_0$ and for all $i \geq 1,\;  G_{i, j} =0 $ if $j \not \in \{ i-1, i+1 \}$, $ G_{i, i-1}= 1-p_i$ and $ G_{i, i+1}= p_i$.

$$ G_{\bar{p}}=
\tiny{
\left[
\begin{array}{cccccccccc}
\!\!1-p_0           \!\!&\!\! p_0 \!\!&\!\!0         \!\!&\!\!0    \!\!&\!\!0            \!\!&\!\!0    \!\!&\!\!0         \!\!&\!\!0    \!\!&\!\!0           \!\!&\!\! \cdots \!\! \\
\!\!1-p_1      \!\!&\!\!0\!\!&\!\!p_1    \!\!&\!\!0    \!\!&\!\!0            \!\!&\!\!0    \!\!&\!\!0         \!\!&\!\!0    \!\!&\!\!0           \!\!&\!\! \cdots \!\! \\
\!\!0               \!\!&\!\!1-p_2    \!\!&\!\!0     \!\!&\!\!p_2  \!\!&\!\!0            \!\!&\!\!0    \!\!&\!\!0         \!\!&\!\!0    \!\!&\!\!0           \!\!&\!\! \cdots \!\! \\
\!\!0   \!\!&\!\!0    \!\!&\!\!1-p_3\!\!&\!\!0\!\!&\!\!p_3  \!\!&\!\!0    \!\!&\!\!0         \!\!&\!\!0    \!\!&\!\!0           \!\!&\!\! \cdots \!\! \\
\!\!\vdots          \!\!&\!\!\vdots \!\!&\!\!\vdots    \!\!&\!\!\vdots \!\!&\!\!\vdots       \!\!&\!\!\vdots \!\!&\!\!\vdots    \!\!&\!\!\vdots \!\!&\!\!\vdots      \!\!&\!\!\ddots
\end{array}
\right]}
$$

It is known (see chapter 5 in \cite{ross}) that $G_p$ is transient if and only if
$$
S_1=
\sum _{n=1}^{+\infty}\frac{(1-p_1)(1-p_2)...(1-p_n)}{p_1 p_2... p_{n}} < \infty \, ,
$$
and positive recurrent if and only if
$$
S_2= \sum_{n=1}^{+\infty} \frac{p_0 p_1... p_{n-1}}{(1-p_1)(1-p_2)...(1-p_n)} < \infty \, .
$$
Thus if both series $S_1$ and $S_2$ do not converge, then $G$ is null recurrent.

Now, consider the sequence
$$
w_n= \frac{(1-p_1)(1-p_3)...(1-p_{n-1})}{p_1 p_3... p_{n-1}} \textrm{ for } n \textrm{ even},
$$
and
$$
w_n= \frac{(1-p_0)(1-p_2)...(1-p_{n-3})(1-p_{n-1})}{p_0 p_2... p_{n-3}p_{n-1}} \textrm{ for } n \textrm{ odd}.
$$

\medskip

\begin{theorem}
\label{passweight}
The following properties hold:

1. If  $X= c_0$ and $\lim w_n=0$ or $X= l^q,\; q \geq 1$ and $\sum_{n=1}^{+\infty} w_n^q <+\infty$, then $G$  is supercyclic on $X$.

2. Let $X \in \{c_0, l^q,\; q \geq 1\}$ and   assume that there exist $n_0 \in \mathbb{N}$ and $\alpha >0$ such that $p_n \geq \frac{1}{2}+ \alpha$ for all $n \geq n_0$,   then there exists $\delta >1$ such that  $ \lambda G$ is frequently  hypercyclic and chaotic for all $\vert \lambda \vert >\delta$.
\end{theorem}

\begin{remark}
Item 1 in Theorem \ref{passweight} implies that there exist null recurrent random walks which are supercyclic on $c_0$.
\end{remark}

\noindent \textbf{Proof:}
1. Assume that $0$ is an eigenvalue of $G$ associated to an eigenvector   $u= (u_n)_{n \geq 0}$. Then
$$
u_1= \frac{p_{0}-1}{p_{0}}u_0 \mbox { and } u_n = \frac{p_{n-1}-1}{p_{n}}u_{n-2} ,\; \forall n \geq 2.
$$
Thus $u_n = (-1)^{n} w_n u_0$ . If $X= c_0$ then $ 0 \in \sigma_{pt}(T)$ if and only if $\lim w_n=0$ .
 If $X= l^q,\; q \geq 1$, then $ 0 \in \sigma_{pt}(G)$ if and only if $\sum_{n=1}^{+\infty} w_n^q <+\infty$.
 In both cases,  we deduce, exactly as done in the proof of Theorem \ref{passeio}, that $G$ is supercyclic on $X$.

2. Assume that $\sum_{n=1}^{+\infty} w_n <+\infty$,

Indeed, let $ v= (v_{i})_{i \geq 0} \in X$ and $u= (u_{i})_{i \geq 0} \in l^{\infty}$ such that $G(u)= v$ and $u_0=0$, then
$$u_1= \frac{v_{0}}{p_0} \quad
\textrm{ and } \quad u_n = \frac{v_{n-1}}{p_{n-1}}+ \frac{p_{n-1}-1}{p_{n-1}} u_{n-2} \textrm{ for all integer } n \geq 2.
$$
Putting $r_n = \frac{p_n -1}{p_{n}}$ for all integer $n \geq 0$,
we obtain by induction that for all integer $n \geq 2$,
\begin{eqnarray*}
u_{n} & = &\frac{1}{p_{n-1}}  v_{n-1}+  \frac{1} {p_{n-3}} r_{n-1}v_{n-3} + \frac{ 1}{p_{n-5}} r_{n-1}r_{n-3} v_{n-5}+ \cdots \\
& & \frac{ 1}{p_{n-2k+1}} ( r_{n-1}r_{n-3} \cdots r_{n-2k+1} )  v_{n-2k-1} + \cdots + \frac{ 1}{p_{t}} ( r_{n-1}r_{n-3} \cdots r_{t+2} ) v_t ,
\end{eqnarray*}
where $t= 0$ if $n $ is odd and $t=1$ otherwise.
Put $S(v)= (0, u_1, \ldots)$ for all $v \in X$, then  $G(Sv)= v$.
Since $p_n \geq \frac{1}{2}+ \alpha$ for all $n \geq n_0$, we have that $r_n \leq \frac{1/2 - \alpha}{1/2+ \alpha} <1$.
We deduce exactly as done in the proof of Theorem \ref{passeio} that there exists $\delta >1$ such that  $ \lambda T$ is frequently  hypercyclic and chaotic for all $\vert \lambda \vert >\delta$. $\square$

\vspace{1em}

{\bf Questions:} 1.  If $X \in \{c_0, l^q,\; q \geq 1\}$ and   $\sum_{n=1}^{+\infty} w_n <+\infty$, can we prove that there exists $\delta >1$ such that  $ \lambda G$ is frequently  hypercyclic and chaotic for all $\vert \lambda \vert >\delta$?

2. If $\sum_{n=1}^{+\infty}  w_n  <+\infty$, then by H\"older inequality, we deduce that
$$\sum_{n=1}^{+\infty} \frac{(1-p_1)(1-p_2) \ldots (1-p_{n})}{p_0 p_1 \ldots p_{n-1}} < +\infty$$ and hence $G$ is transient.

Does there exist $G$ transient and not supercyclic on $l^1$ such that $\sum_{n=1}^{+\infty} w_n <+\infty$?

\medskip

\begin{theorem} \label{notGRW}
If $X = c_0$  and $\sum_{n=1}^{+\infty}  w_n ^{-1} <+\infty$ or  $X= l^1$ and $1/w_n$ is bounded or $X= l^q,\; q > 1$ and  $\sum_{n=1}^{+\infty} ( w_n )^{-\frac{q}{ q-1}} <+\infty$,  then  $\lambda G$ is not hypercyclic for all $\vert \lambda \vert >1$.
\end{theorem}

\begin{remark}
Theorem \ref{notGRW} is the closest we get to Theorem \ref{notSRW}. We conjecture that $G$ is not supercyclic under the hypothesis of Theorem \ref{notGRW}.
\end{remark}

\noindent \textbf{Proof:}
Assume that  $ 0 $ is an element of $\sigma_{pt}(X', G')$ and $u= (u_n)_{n \geq 0}$ an eigenvector associated to $0$, then $u T= 0$.
Thus $(1-p_0)u_0 + (1- p_1)u_1= 0$ and $p_n u_n + (1- p_{n+2}) u_{n+2}= 0$ for all $n \geq 0$.
Hence for all $n \geq 1$, we have
 $$u_{2n} = \frac{p_{2n-2} p_{2n-4} \ldots p_{0}}{(p_{2n-2}-1)(p_{2n-4}-1) \ldots (p_{0}-1)} u_0$$ and
$$u_{2n+1} = \frac{p_{2n-1} p_{2n-3} \ldots p_{1}}{(p_{2n-1}-1)(p_{2n-3}-3) \ldots (p_{1}-1)} \frac{(1-p_{0})}{p_1 -1} u_0.$$
Hence if  $X = c_0$  and $\sum_{n=1}^{+\infty}  w_n ^{-1} <+\infty$ or  $X= l^1$ and $1/w_n$ is bounded or $X= l^q,\; q > 1$ and  $\sum_{n=1}^{+\infty} ( w_n )^{-\frac{q}{1- q}} <+\infty$, we have  $0 \in \sigma_{pt}((\lambda G)^{'}, X'))$.
Thus, by Lemma \ref{spectrhyper}, $\lambda G$ is not hypercyclic  for all $\lambda$. $\square$

\medskip

{\bf Question:}
If $G$ is null  recurrent and $X= l^1$, Is  $\lambda T$ is not hypercyclic for all $\vert \lambda \vert >1$?

\medskip

\begin{remark}

Let  $\bar{p} = (p_{n})_{n \in \mathbb{Z}}$ be a sequence of probabilities, and denote by $\overline{G}_{\bar{p}} := \overline{G}$ the transition operator of the spatially inhomogeneous simple random walks on $\mathbb{Z}$, defined by: For all $i \in \mathbb{Z},\; \overline{G}_{i, i-1}= 1-p_i,\; \overline{G}_{i, i+1}= p_i$ and $\overline{G}_{i, j}= 0$ if $j \not \in \{i-1, i+1\}$.
Then by using the same method done in Theorem \ref{passweight}, we can prove the following results:

1. If $\lim  w_n=0$ and   $\lim w_{-n}^{-1}=0$ then $\overline{G}$  is supercyclic on $c_0$.
If $q \geq 1$ and $\sum_{n =1}^{+\infty} w_n^{q} <+\infty$ and $ \sum_{n =1}^{+\infty} w_ {-n}^{-q} <+\infty$, then $\overline{G}$  is supercyclic on $l^q$.

2. Let $X \in \{c_0, l^q,\; q \geq 1\}$ and   assume that there exist positive constants $n_0, n_1, \alpha $ such that $p_n \leq \frac{1}{2}- \alpha$ for all $n \geq n_0$ and
 $p_n \geq \frac{1}{2}+ \alpha$ for all $n \leq -n_1$,   then there exists $\delta >1$ such that  $ \lambda \overline{G}$ is frequently  hypercyclic and chaotic for all $\vert \lambda \vert >\delta$.
\end{remark}

\textbf{Question:} Does there exist a transient Markov operator that is not supercyclic?
Does there exist a null recurrent Markov operator that is supercyclic on $l^1$?

\medskip

\noindent{\bf Acknowledgment:} The authors would like to thank El Houcein El Abdalaoui and Patricia Cirilo  for fruitful discussions.


\begin{thebibliography}{99}

\bibitem{acmv} H. Elabadaloui, P. Ciriclo, A. Messaoudi, G. Valle, {Dynamics of stochastic adding machines}, in preparation.

\bibitem{bm} F. Bayart, E. Matheron, {\it Dynamics of linear opertors}, Cambridge university Press, 2009.

\bibitem{BerBonMulPer13}
    N. C. Bernardes Jr., A. Bonilla, V. M\"uller and A. Peris,
    {\it Distributional chaos for linear operators},
    J. Funct.\ Anal.\ {\bf 265} (2013), no.\ 9, 2143--2163.

\bibitem{BerBonMulPer14}
    N. C. Bernardes Jr., A. Bonilla, V. M\"uller and A. Peris,
    {\it Li-Yorke chaos in linear dynamics},
    Ergodic Theory Dynam.\ Systems {\bf 35} (2015), no.\ 6, 1723--1745.

\bibitem{BesMenPerPui16}
    J. B\`es, Q. Menet, A. Peris and Y. Puig,
    {\it Recurrence properties of hypercyclic operators},
    Math.\ Ann.\ {\bf 366} (2016), no.\ 1, 545-572.

\bibitem{cg} Carleson, L.; Gamelin, T.: {\it Complex Dynamics}, Springer, 1993.

\bibitem{em} K.G. Erdmann, A.P. Manguillot, {\it Linear Caos}, Springer.

\bibitem{k} B. Kitchens, {\it Symbolic Dynamics}, Springer, 1997.

\bibitem{ross} S. Ross, {\it Stochastic Processes,} Second edition, Wiley \& Sons, 1996.

\bibitem{r} W.~Rudin, {\it Functional analysis,} Second edition. International Series in Pure
and Applied Mathematics. McGraw-Hill, Inc., New York, 1991.

\bibitem{SGriEMat14}
    S. Grivaux and \'E. Matheron,
    {\it Invariant measures for frequently hypercyclic operators},
    Adv.\ Math.\ {\bf 265} (2014), 371--427.

\bibitem{SGriMRog14}
    S. Grivaux and M. Roginskaya,
    {\it A general approach to Read's type constructions of operators
    without non-trivial invariant closed subspaces},
    Proc.\ Lond.\ Math.\ Soc.\ (3) {\bf 109} (2014), no.\ 3, 596--652.

\bibitem{Y}
    K. Yosida,
    {\it Functional Analysis}, Sixth Edition,
    Springer-Verlag, Berlin and New York, 1980.

\end{thebibliography}
\end{document}